# Some Statistical Tools to Measure The Effectiveness of a (LEM) Model on Creative Thinking Among Intermediate Students in Mathematics

Hiba Ali Kareem[1]    Abbas N. A. Ameer[2]    Maan A. Rasheed[3*]

[1,2,3] Department of Mathematics, College of Basic Education, Mustansiriyah University, Baghdad, Iraq

*maan.rasheed.edbs@uomustansiriyah.edu.iq

**Abstract**

The article aims to identify the effectiveness of the (LEM) model for a multimodal environment on creative thinking among first grade intermediate students in Mathematics using some statistical tools. We achieved the credibility, stability and prepared a test of creative thinking for its three levels (fluency, flexibility and originality). We achieved their credibility by presenting the results to a group of experts and arbitrators. The stability (92%) was extracted using Alpha Cronbach formula. After applying the test and processing the data statistically using t-test for the independent sample, the results were as follows: There was a statistically significant difference at the level of significance (0.05) between the average scores of the experimental and control group students and for the experimental group in the test of creative thinking. We conclude that the adoption of the (LEM) model for a multimodal environment has a positive effect on creative thinking in mathematics.

## 1. INTRODUCTION

The educational environment is one of the basic pillars in the educational process. It must adopt all the future hopes to improve the educational process. As far as the level of the environment is concerned, as far as the growth and development of the learners is concerned. The educational environment is the one that inspires activity inside and outside the classroom and adds to the book, content, activities and means that complement the shortage in order to serve the learner. Most educators agree that education for thinking or learning skills is an important goal of education, and schools must do whatever they can to provide learners with the opportunity to think.

Our schools rarely provide learners with the opportunity to carry out educational tasks stemming from their curiosity or based on questions raised by them. The researchers noted that most of the staff in the educational field emphasize that the mission of the school is not to fill the minds of the students with information as much as to stimulate thinking and creativity. They coexist with the prevailing practices in our schools and none of them tried to break the familiar wall or get out of it.

In addition, some educational bodies and mathematics supervisors (after meeting with the researchers) stressed the need for learners in the middle stage to develop a clean and mature educational environment that promotes the development of thinking, especially creative thinking. Among these models, which may work and contribute to providing an educational environment that facilitates learners' learning, learning, and self-knowledge as well as their daily use is the LEM (Learning Environment Modeling) model, for a multimodal environment, which may reflect positively on learners' achievement and creative thinking. So that, we see many learners do not think well not because they lack intelligence, but because they did not find the right direction and the educational environment necessary for them and that shows the inability of the methods of teaching used to achieve the goal of thinking or development of learners in the daily ration.

Through the our knowledge of the literature, theses and articles that dealt with creating a model educational environment and the recent changes that are sweeping the world in recent years, we noticed that many countries have re-examined their educational systems in general and the emphasis on the educational environment in particular, and provide them with professional skills and modern educational environments, commensurate with the nature of the learner, [1].

In this article we use some statistical tools in order to identify the effectiveness of the (LEM) model for a multimodal environment on creative thinking among first grade intermediate students in Mathematics.

To achieve the aims of the research, the following null hypothesis has been set:
There is no statistically significant difference at (0.05) level of significance between the mean of scores of the experimental group which were taught mathematics according to (LEM) Model and that of the control group which was taught by the adopted method in the creative thinking test.

$$H_O: \bar{x}_1 = \bar{x}_2$$
$$H_1: \bar{x}_1 \neq \bar{x}_2$$

In order to verify the validity of these null hypotheses, the researchers conducted their experiment, where the study was limited to the first-grade intermediate students. The experimental design was adopted with partial control of two groups (experimental and control) with the post-test. We have chosen (Al-Hauraa intermediate school for girls, Baghdad-Iraq) to apply the experiment. Class (A) was chosen randomly to represent the experimental group, while Class (e) represented the control group. The sample of the research consists of (75) students divided to (35) students in the experimental group and (40) students in the control group.

We have equalized between the students of the two research groups statistically using the TI test for two independent samples in the following variables: (intelligence test, previous achievement in mathematics, previous knowledge).

After we identified the subjects of the scientific material to be studied during the duration of the experiment, the researchers formulated the behavioral goals, and prepared the teaching plans, and presented to a group of experts and specialists to judge the validity, and made the necessary adjustments and the plans are ready to be applied in the light of their views.

we achieved its credibility and stability and prepared a test of creative thinking for its three levels (fluency, flexibility and originality) component (20) paragraph. We achieve the credibility by presenting the results to a group of experts and arbitrators. The stability (92%) was extracted using Alpha Cronbach formula.

After applying the test and processing the data statistically using t-test for the independent sample, the results were as follows:
There was a statistically significant difference at the level of significance (0.05) between the average scores of the experimental and control group students and for the experimental group in the test of creative thinking and collection.

## 2. DEFINITIONS AND CONCEPTS

**2.1 Effectiveness**: In [2], Effectiveness was defined as: measuring the level of achievement or achievement of goals.

**2.2 Learning Environments Modeling (LEM)**

(LEM) is a way to help teachers enhance their understanding of the learning environments they create, improving design decisions and providing clear and effective communication [3]. It offers a new way to design a course that transfers the design experience.

It is not a rigid or specific process, but it is a flexible, changeable idea and an educational design platform. It is also a tool to understand how to design learning environments, to make deliberate design decisions, and to share those decisions with others involved in curriculum design.

LEM uses new and effective ways to represent how learning environments are designed in a similar way to how the architect designs physical spaces.

**2.3 Learning Environments Modeling Language (LEML)**

LEML equips the learning environment architect with a language and visual toolkit for communicating complex information about how learning environments are designed. LEML is a language and framework that consists of a set of symbols that can be assembled together to represent information about the design of learning environments. LEML reframes the design experience by providing a concise and easy-to-understand design language for learning environments. Leveraging this user-friendly design language, LEML helps learning environment architects enhance communication and improve the value of learning environment design experiences.

LEML is comprised of four primary features: Building Blocks, Contexts, Actions, and Notations.

- Building blocks represent the core components, or system nodes, of a learning environment. LEML has five types of building blocks that can be configured to represent the design of any learning environment. These are: Information, Dialogue, Feedback, Practice, and Evidence.

- The Contexts in LEML describe the mode, or space in which the building blocks reside. There are four contexts represented in LEML: Classroom, Online Synchronous, Online Asynchronous, and Experiential. The contexts are represented by color-coded boxes.

- The Actions in LEML identify the connections and transitions between building blocks and who/what is responsible for performing those actions .There are three types of actions in LEML :Learner Action, Facilitator Action, and System Action.

- Notations can be added to enhance the meaning and usefulness of the learning environment model. While there are numerous possibilities for useful notations, there are two standard notation elements used: Start-Stop and Objective ID.

**2.4. The Multimodal Environment Model**

This model describes the integration of four different types of learning methods into a single learning process [4].

In this model the class begins with a written overview on the subject, as well as some videos are viewed on the Internet that offer a concept of the subject and show skills. Students then meet in the classroom to work on (question groups) and receive feedback from the teacher on their work in class. After that, students meet in the lab to observe skills in practice. Learners then participate in online discussion of ideas (synchronously or asynchronously) and each student presents individual solutions to the issue via e-mail or Dropbox. The teacher must keep working hours online through the electronic part of the lesson.

The Figure (1) illustrates how to integrate, into this model, four different types of learning frameworks in a single learning process

We have found in this model the possibility of application in the classroom in Iraqi schools. Therefore, we have chosen this model for the subject of research to see its effectiveness in the creative thinking of students in the first grade intermediate in mathematics.

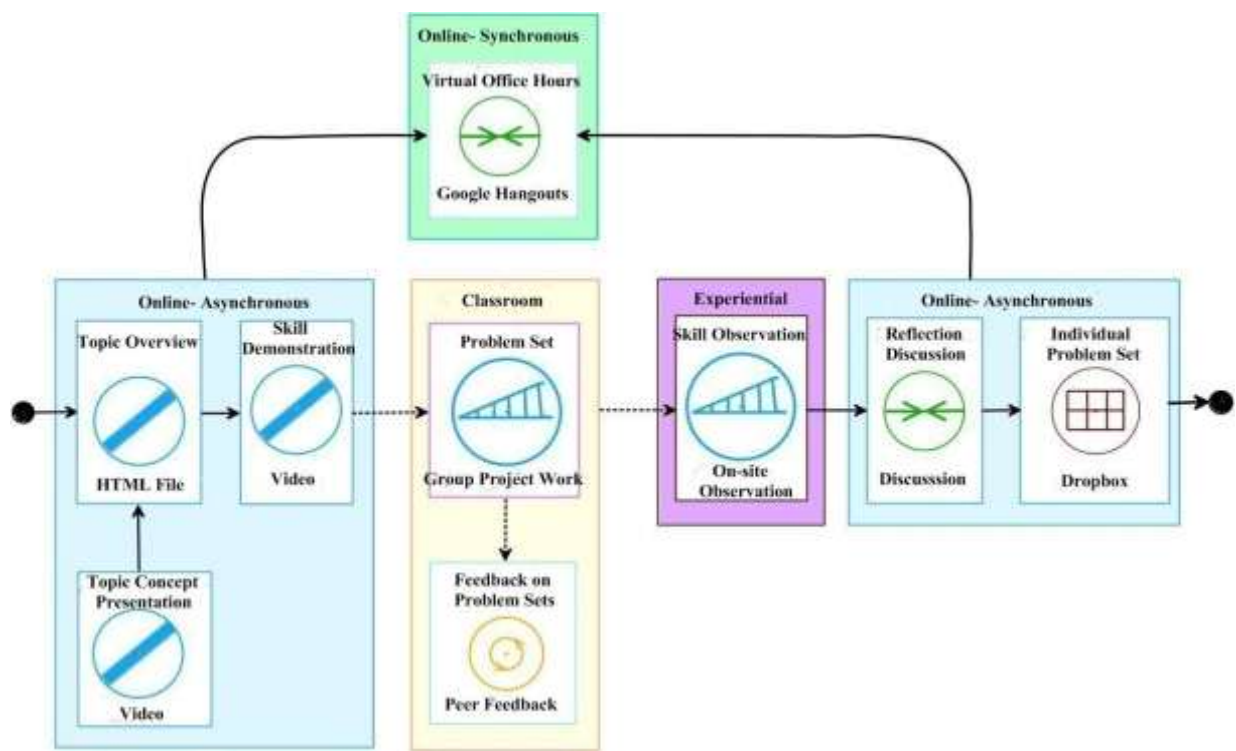

**Figure 1: Multimodal Environment Model**

## 2.5 Creative thinking

Creative thinking is a mental activity, is an exception that stems from a problem or situation that is attractive to attention, a constant that moves its owner from one location to another and from one solution to the other without the need to walk in a routine way, change is his style and purpose, [5].
"Creating or producing new ideas or finding new solutions to challenges," Zubaidi said, [6].

### 2.5.1 Creative Thinking Skills:
Through the study of educational literature and research, We have noted that creative thinking has several skills and that the skills that suit the students who has been studied, these skills are fluency, flexibility and originality, so the researchers have dealt with these three skills. The skills of creative thinking are defined as: "A new purposeful and goal-oriented production, which is the ability of the mind to form new relationships that change the reality of the individual, where conservation and memorization go beyond thinking, study, analysis, conclusion, innovation and creativity." ,[7]. Creative individuals have many creative abilities to make creative production. Many studies and research have revealed the most creative abilities that determine the creative potential of individuals: fluency, flexibility, and originality.
We have tried to take the skills of creative thinking that can be developed or developed during the current research. After presenting the skills of creative thinking to a group of specialized arbitrators, a test of creative thinking was agreed upon which deals with skills (fluency, flexibility, originality) to identify the effect of using a renewable and multimodal learning environment in these skills.

## 3. STATISTICAL TOOLS

After the data collection and analysis using the SPSS to measure the variables of research in the students of the research sample (first grade students average), We resorted to the use of statistical tools in the analysis of research data and the extraction of results, and used the following statistical tools:

A) Difficulty factor for paragraphs, [8]: used to calculate the coefficient of difficulty of paragraphs to test creative thinking in Mathematics:

$$D.C = \frac{P_u - P_I}{x(n_u + n_I)} \qquad (1)$$

where

$P_u$ : Total scores obtained by students of the upper group.
$P_I$ : Total scores obtained by the lowest group students.
$x$: The largest number of grades obtained by students in question.
$n_u$: Number of students in the upper group who answered the question correctly.
$n_I$: The number of students in the lower group who answered the question is correct.

B) The Alpha-Cronbach equation, [9]: Used to calculate the stability of the test of creative thinking in Mathematics.

$$\alpha = \frac{n}{n-1}\left[1 - \frac{\sum s_i^2}{s_x^2}\right], \qquad (2)$$

where
$\alpha$: stability coefficient
$N$: Number of paragraphs,
$\sum SI^2$: Variation of scores on each paragraph in the test,
$SX^2$: Variation of scores on the test as a whole.

C) Pearson correlation coefficient, [10]: Use to extract the correlative relationship between the exploratory sample for each of the paragraphs of the creative thinking test with the total score for testing and extraction (stability by retesting).

$$R = \frac{n\sum xy - (\sum x)(\sum y)}{\sqrt{[n\sum x^2 - (\sum x)^2][n\sum y^2 - (\sum y)^2]}} \qquad (3)$$

where
R.: correlation coefficient
n: Student number
x,y: Variables

D) Cooper's equation, [11]: Adopted to calculate the proportion of agreement between the arbitrators.

$$P = \frac{NP}{NP+NNP} \qquad (4)$$

where
P: ratio of agreement
NP: number of times agreed
NNP: The number of times it is not agreed

E) The equation of determining the size of the effect (Claes equation), (Allam,1989):
We used the Claas equation to demonstrate the effect of the independent variable (LEM) for the multimodal environment in the dependent variable (creative thinking).

$$E.K = \frac{\overline{x_1}-\overline{x_2}}{S_2} \qquad (5)$$

where
$E.K$ : Claes equation
$\overline{x_1}$ : The arithmetic mean of the experimental group.
$\overline{x_2}$ : Arithmetic mean of the control group.
$S_2$ : Standard deviation of the control group.

## 4. RESULTS AND DISCUSSION

The next tables 1-3 show the statistical results obtained by using the statistical tools eqs. 1-5 mentioned above with SPSS software.

**Table 1**
**Difficulty, ease, and discrimination for test paragraphs to measure creative thinking in Mathematics**

| NO. | Sum. of Marks | | Difficulty coefficient | Discriminatory power |
|---|---|---|---|---|
| | Max. | Min. | | |
| 1 | 120 | 62 | 0.67 | 0.42 |
| 2 | 113 | 66 | 0.66 | 0.34 |
| 3 | 111 | 53 | 0.60 | 0.42 |
| 4 | 105 | 46 | 0.56 | 0.43 |
| 5 | 112 | 63 | 0.64 | 0.36 |
| 6 | 110 | 47 | 0.58 | 0.46 |
| 7 | 133 | 65 | 0.73 | 0.50 |
| 8 | 125 | 81 | 0.76 | 0.32 |
| 9 | 129 | 80 | 0.77 | 0.36 |
| 10 | 112 | 59 | 0.63 | 0.39 |
| 11 | 101 | 53 | 0.57 | 0.35 |
| 12 | 133 | 73 | 0.76 | 0.44 |
| 13 | 124 | 81 | 0.75 | 0.31 |
| 14 | 105 | 63 | 0.62 | 0.31 |
| 15 | 100 | 45 | 0.53 | 0.40 |
| 16 | 109 | 67 | 0.65 | 0.31 |
| 17 | 89 | 33 | 0.45 | 0.41 |

| 18 | 114 | 69 | 0.67 | 0.33 |
| 19 | 91 | 41 | 0.48 | 0.37 |
| 20 | 117 | 64 | 0.67 | 0.39 |

**Table (2)**
**The correlation coefficient values for each paragraph and its field and the total number of creative thinking selection**

| No. | correlation coefficient | No. | correlation coefficient | Level of significance |
|---|---|---|---|---|
| 1 | 0.344 | 11 | 0.636 | 5% |
| 2 | 0.522 | 12 | 0.526 | 5% |
| 3 | 0.757 | 13 | 0.430 | 5% |
| 4 | 0.622 | 14 | 0.684 | 5% |
| 5 | 0.631 | 15 | 0.642 | 5% |
| 6 | 0.660 | 16 | 0.653 | 5% |
| 7 | 0.655 | 17 | 0.697 | 5% |
| 8 | 0.759 | 18 | 0.707 | 5% |
| 9 | 0.679 | 19 | 0.814 | 5% |
| 10 | 0.619 | 20 | 0.717 | 5% |

**Table (3)**
**Test results (t-test) to determine the difference between the mean scores of the two research groups (Experimentation and control) in the test of creative thinking as a whole**

| Group | No. | Standard Mean | Standard Deviation | Values (t-test) | | Freedom degree | ETA square (impact size) | Level of significance |
|---|---|---|---|---|---|---|---|---|
| | | | | Calculated | Tabular | | | |
| **Experimentation** | 35 | 77.8286 | 17.70591 | 3.564 | 1.99 | 73 | 0.148 | 0.05 |
| **control** | 40 | 63.0250 | 18.15106 | | | | | |

To verify the null hypothesis, there was no statistically significant difference at the level of (0.05) between the average score of the students of the experimental group studied according to the LEM model for the multimodal environment and the average score of the control group students studied according to (the usual method) in the creative thinking test students of the experimental and control experimental groups and the mean differences of test scores for the experimental and control groups in the test of creative thinking as a whole.

The mathematical averages and standard deviations of female students 'scores in the creative thinking test and the three abilities as a whole (fluency, flexibility, originality) of the experimental group were derived and then compared to the mean and standard deviations of the students' grades in the creative thinking test of the control group. The results showed that the experimental averages of the experimental group were higher than the averages Arithmetic of the control group. It is clear from Table (3) that the average level of creative thinking in the experimental group (23.6857) and its standard deviation (5.37845) (T-test) for two separate independent scores. The calculated T value (3.564) is greater than the scale value of (1.99 at the level of significance (0.05) and the degree of freedom (73). This means that there is a statistically significant difference at the level of (0,05) in the positive thinking of the students who studied according to the LEM model for the multimodal environment compared to the thinking of the students who studied in the usual way for the experimental group in the three capacities. The null hypothesis is rejected, and the value of the ETA box is 0.151. Thus, the students of the research sample in general have the three abilities (fluency, flexibility, originality) of creative thinking in mathematics at the required level.

The results of this study can explain the results of this variable, which showed the superiority of the students of the experimental group who studied according to the LEM model of the multimodal environment to the students of the control group who studied according to the usual method of teaching in the creative thinking test.

## 5. CONCLUSIONS

In light of the results of the research, we concluded that the adoption of the LEM model for a multimodal environment has a positive effect on creative thinking in mathematics. Therefore, the mathematics teachers should adopt the (LEM) model for a multimodal environment in teaching because of its positive impact on the development of creative thinking and encourage teachers to pay attention to teaching thinking as a mental activity that helps to move learning into practice.

**ACKNOWLEDGEMENTS**

The authors would like to thank Mustansiriyah University (www.uomustansiriyah.edu.iq) Baghdad-Iraq for its support in the present work.

**REFERENCES**

1. Saada Wadat Ahmed. (2015): Thinking and Learning Skills, Dar Al Masirah Publishing, Distribution and Printing, Amman.
2. Bruce A. Kirchhoff, (1977), Organization Effectiveness Measurement and Policy Research, The Academy of Management Review, vol. 2, No. 3, pp. 347-355.
3. Bucky Dodd : John Gillmore, (2018): Learning Environment Modeling, Institute for learning environment design, university of central Oklahoma.
4. Bucky Dodd : John Gillmore, (2017): The learning Designer's: Guide to LEM, Institute for learning environment design, university of central Oklahoma.
5. Al-Khalili, Amal (2005): Developing Children's Innovation Capabilities, i 1, Dar Safaa Publishing House, Amman.
6. Zubaidi, Khoula (2006): Thinking Skills and Problem Solving Skills; Al-Shukri Library, Riyadh.
7. Ghadana Saeed Al - Binali, (2003), The extent to which teachers used social studies for thinking skills in teaching primary school students in Qatar, Gulf Message Magazine, Riyadh, Issue 99, pp 69-111.
8. Allam, Salahuddin Mahmoud (2007): Measurement and Evaluation in the Teaching Process, 1, Dar Al-Masirah, Amman.
9. Al-Nabhan, Mousa (2004): The Basics of Measurement in Behavioral Sciences, I 1, Dar Al-Shorouk for Publishing and Distribution, Amman.
10. Al - Bayati, Abdul - Jabbar Tawfiq (2008): Statistics and its applications in educational and psychological sciences, $1^s$ edition, Dar Athra for Publishing, Distribution and Printing, Amman.

11. Cooper Hohn charles, (1974): Measurement and Analysis of Behavioral Teachingues chio, Emeirll Columbus.
12. Allam, Salah Eddin Mahmoud (1989): Design and experimentation of the model of the educational system for the competencies of psychological statistics using the entrance of the calendar of referees Reference (Journal of Social Sciences) Vol. (7), No. (3), Al-Azhar University.